\documentstyle{amsart}
\newcommand{\N}{\text{\normalshape\sf N}}
\newcommand{\B}{\text{\normalshape\sf B}}
\newcommand{\Bzero}{\text{\normalshape\sf B$_0$}}
\newcommand{\A}{\text{\normalshape\sf A}}
\newcommand{\D}{\text{\normalshape\sf D}}
\newcommand{\pD}{\text{\normalshape\sf {D}$_\sigma$}}
\renewcommand{\H}{{\text{\normalshape\sf H}}}

\newcommand{\lbv}{[\![}
\newcommand{\rbv}{]\!]}
\newcommand{\reals}{{\Bbb R}}
\newcommand{\rationals}{{\Bbb Q}}

\newcommand{\dom}{{\text{\normalshape\sf {dom}}}}

\newcommand{\QED}{\hspace{0.1in} \Box \vspace{0.1in}}

\newcommand{\<}{\langle}
\renewcommand{\>}{\rangle}

\newtheorem{theorem}{Theorem}
\newtheorem{lemma}[theorem]{Lemma}
\newtheorem{corollary}[theorem]{Corollary}
\newtheorem{definition}[theorem]{Definition}

\theoremstyle{definition}

\newcommand{\lesdot}{\mathrel{\mathord{<}\!\!\raise 
0.8 pt\hbox{$\scriptstyle\circ$}}}
\newcommand{\Proof}{{\sc Proof} \hspace{0.2in}}

\newcommand{\lft}[2]{\mathopen\ifcase#1{}\oo\or
                        \big#2\or\Big#2\else\oo\fi} 
\newcommand{\rgt}[2]{\mathclose\ifcase#1{}\oo\or
                        \big#2\or\Big#2\else\oo\fi} 
\begin{document}

\title{Remarks on sets related to trigonometric series}
\author{Tomek Bartoszy\'{n}ski}
\address{Department of Mathematics\\
   Boise State University\\
   Boise, Idaho 83725}
\email{tomek@@math.idbsu.edu}
\author{Marion Scheepers}
\address{Department of Mathematics\\
   Boise State University\\
   Boise, Idaho 83725}
\email{marion@@math.idbsu.edu}
\thanks{The second author was partially supported by Idaho State Board
  of Education Grant 94-051.}
\thanks{The authors would like to thank L. Bukovsk\'y and M. Repick\'y for their helpful
  comments.}
\maketitle
\begin{abstract}
   We show that several classes of sets, 
   like $\N_0$-sets, Arbault sets, $\N$-sets and pseudo-Dirichlet sets are 
   closed under adding sets of small size. 
\end{abstract}
\section*{Introduction}
The goal of this paper is to prove the theorem below. In this first
section we will introduce 
all necessary definitions. For more information we refer the reader to
the survey paper \cite{Buki}.

\begin{theorem}\label{first}
   Let $X$ be a set of real numbers. 
\begin{enumerate}
\item If $A$ is an $\N$-set and if $|X|< \frak t$, then $A \cup X$ and
  $A+X$ are 
  $\N$
sets. 
\item If $A$ is an $\N_0$-set and if $|X|< \frak h$, then $A \cup X$ and
  $A+X$ are
  $\N_0$ sets. 
\item If $A$ is a $\pD$ set and if $|X| < \frak h$, then 
   $A \cup X$ and
  $A+X$ are  $\pD$ sets.
\item If $A$ is an $\A$ set and if $|X| < \frak s$, then $A \cup X$
  and
  $A+X$ are  $\A$ sets. 
\item If $A$ is an $\H_{\sigma}$ set and if $|X|<\frak s$, 
then $A\cup X$ is an
  $\H_{\sigma}$-set.  
\end{enumerate}
\end{theorem}

Note that part $(1)$ strengthens the main result of \cite{BuBu},
part (2) generalizes the result of Arbault and Erd\"{o}s and (4)
strengthens the result from \cite{Kch}.

Parts (1)-(4) of the following definition are classical. For more
information see 
\cite{bary64}.  Parts (5) and (6) were introduced in  \cite{elias93}.
\begin{definition}\label{defdef}
Suppose that $X \subseteq [0,1]$.
\begin{enumerate}
\item $X$ is an $\N$ set if there exists a sequence $\<a_n: n \in
  \omega\>$ such that 
the series $\sum_{n=0}^\infty a_n \cdot \sin(\pi \cdot n \cdot x)$
converges absolutely 
for all $x \in X$ and $\sum_{n=0}^\infty a_n=\infty$.
\item $X$ is an $\N_0$ set if there exists an increasing function $f
  \in \omega^\omega$ 
such that the series 
$\sum_{n=0}^\infty \sin\lft1(\pi \cdot f(n) \cdot x\rgt1)$ 
converges absolutely
for all $x \in X$.
\item $X$ is an $\A$ (Arbault) set if there exists an increasing
  function $f \in \omega^\omega$ 
such that $\lim_{n \rightarrow \infty} \sin\lft1(\pi \cdot f(n) \cdot
x\rgt1) = 0$ for every $x \in X$.
\item $X$ is a $\D$ (Dirichlet)  set if there exists an increasing function 
$f \in \omega^\omega$
such that the sequence
$\<\sin\lft1(\pi \cdot f(n) \cdot x\rgt1) : n \in \omega\>$ 
converges uniformly to $0$ on $X$.
\item $X$ is a $\pD$ (pseudo Dirichlet)  set if $X$ is a countable
  union of $\D$ sets.
\item $X $ is a $\Bzero$ set if there exists an increasing function 
$f \in \omega^\omega$ and a number $n_0 \in \omega$ such that for
every $x \in X$,
$\sum_{n \geq n_0} |\sin\lft1(\pi \cdot f(n) \cdot x\rgt1)| \leq n_0$.
\item $X $ is a $\B$ set if there exists a sequence of positive reals 
$\<a_n: n \in
  \omega\>$ with $\sum_{n=0}^\infty a_n=\infty$ and there is a
  number $n_0$ such that for all $x \in X$,
$\sum_{n > n_0}^\infty a_n \cdot |\sin(\pi \cdot n \cdot x)| \leq n_0$.
\end{enumerate}
\end{definition}
It is easy to see that 
\begin{lemma}
$X$ is a $\pD$  set iff there exists a  sequence $\<\varepsilon_n : n
\in \omega\>$ converging to zero and an increasing function
$f \in \omega^\omega$ such that 
$$\forall x \in X \ \forall^\infty n \ 
\lft2|\sin\lft1(\pi \cdot f(n) \cdot x\rgt1)\rgt2| <\varepsilon_n.~\QED$$
\end{lemma}

For $x \in \reals$ let $\lbv x \rbv$ be 
the distance from $x$ to the nearest integer.
Note that $\lbv x \rbv \leq  |\sin(\pi \cdot x)| \leq \pi \cdot \lbv x \rbv$.
Therefore we can replace every occurrence of the function $\sin(\cdot)$ by 
$\lbv \cdot \rbv$ in the definitions above.

\begin{definition}
A set $X \subseteq [0,1]$ is an $\H$ (Hardy) set if
there exist  reals $a \in (0,1)$, $\varepsilon < 1/2$ and an
increasing function 
$f \in \omega^\omega$ such that
$$\forall x \in X \ \forall n \ \lbv f(n) \cdot x - a\rbv < \varepsilon.$$
$X$ is an $\H_\sigma$ set if $X$ is a countable union of $\H$ sets.
\end{definition}
The terminology $\H$-set is due to Rajchmann, and is in honour of
   Hardy and Littlewood. It is well-known that an $H$-set has measure zero
   and is of the first category. A translation of an $H$-set is an $H$-set.

\begin{lemma}[Eli\'a\v{s} \cite{elias93}]
We have the following inclusions:
$$\begin{array}{ccccc}
\D & \rightarrow & \pD & \ & \ \\
\downarrow & \ & \downarrow & \ & \ \\
\Bzero & \rightarrow & \N_0 & \rightarrow & \A \\
\downarrow & \ & \downarrow & \ & \downarrow \\
\B & \rightarrow & \N & \rightarrow & P(\reals)
\end{array}$$
Moreover $\A \subset \H_\sigma$.~$\QED$
\end{lemma}

We need  definitions of the following cardinal invariants:
\begin{definition}
For $X,\ Y \subseteq \omega$ let $X \subseteq^\star Y$ denote
 that $X \setminus Y$
is finite.
\begin{enumerate}
\item $\frak p$ is the least size of a family $\cal A \subseteq
  [\omega]^\omega$ 
such that for every $A_1, \ldots, A_n \in \cal A$, $A_1 \cap \cdots
\cap A_n$ is  
infinite and there is no set $B \in [\omega]^\omega$ 
such that $B \subseteq^\star A$ for
all $A \in \cal A$,
\item $\frak t$ is the least size of a family $\{A_\alpha:
  \alpha<\kappa\} \subseteq [\omega]^\omega$
such that $A_\alpha \subseteq^\star A_\beta$ for $\alpha \geq \beta$
and such that there is no set $B \in [\omega]^\omega$ such that 
$B \subseteq^\star A_\alpha$ for
all $\alpha<\kappa$,
\item $\frak h$ is the least size of a family 
$\{\cal A_\alpha: \alpha<\kappa\}$ such that  ${\cal A}_\alpha
\subseteq [\omega]^\omega$ and 
$$ \forall X \in [\omega]^\omega \ \forall \alpha<\kappa \ \exists A
\in \cal A_\alpha \ 
A \subseteq^\star X$$
but there is no set $B \in [\omega]^\omega$ such that 
$$\forall \alpha<\kappa \ \exists A \in \cal A_\alpha \ B \subseteq^\star A.$$
\item $\frak s$ is the least size of a family $\cal A \subseteq
  [\omega]^\omega$ 
such that  for every set $B \in [\omega]^\omega$ there is an $A \in {\cal
  A}$ such 
that $|B \cap A|=\lft1|B \cap (\omega \setminus A)\rgt1|=
\boldsymbol\aleph_0$.
\item $\frak b$ is the least size of a family $F \subseteq
  \omega^\omega$ such that there is no $g \in \omega^\omega$ such that
  $$\forall f \in F \ \forall^\infty n \ f(n)<g(n).$$
\end{enumerate}
\end{definition}

Note that $\frak h$ is the smallest size of a family of open dense
subsets of $P(\omega)/\text{\sf{fin}}$ which has empty intersection.
 
It is well known that $\frak p \leq \frak t \leq \frak h \leq \frak
s$, $\frak h \leq \frak b$ and 
consistently $\frak t < \frak h$, $\frak h < \frak s$ and $\frak h <
\frak b$ (see
\cite{vdo:inttop} and \cite{Vau90Sma}).

The following lemma shows that in Theorem \ref{first} we need only
be concerned with unions of sets.
\begin{lemma}
  Suppose that $A,B$ are two nonempty sets of reals.
  \begin{enumerate}
 \item    $A \cup B$ is an $\N$ set iff $A+B$ is an $\N$ set,
\item    $A \cup B$ is an $\N_0$ set iff $A+B$ is an $\N_0$ set,
\item    $A \cup B$ is a $\pD$ set iff $A+B$ is a $\pD$ set,
\item    $A \cup B$ is an $\A$ set iff $A+B$ is an $\A$ set,
\item    $A \cup B$ is an $\B$ set iff $A+B$ is an $\B$ set,
\item    $A \cup B$ is an $\Bzero$ set iff $A+B$ is an $\Bzero$ set.
  \end{enumerate}
\end{lemma}
\Proof We will only prove part $(1)$. 

Suppose that $A +B$ is an $\N$ set.  $A \cup B$ is contained
in a translation of  $A +B$ that is an $\N$ set.

\vspace{0.1in}

To show that  $A +B$ is an $\N$ set whenever  $A  \cup B$ is an $\N$ set
use the fact that
$$\lft1|\sin(x+y)\rgt1| = \lft1| \sin(x)\cos(y) + \sin(y)\cos(x)\rgt1|
\leq
\lft1|\sin(x)\rgt1| + \lft1|\sin(y)\rgt1|.~\QED$$

Note that by a theorem of Marcinkiewicz (see \cite{bary64}, chapter
12.11), there are 
two $\N_0$ sets $A$ 
and $B$ such that $A+B$ is not an $N_0$ set.

\section{{\normalshape\sf N} sets}
In this section we will prove the first part of Theorem \ref{first}.

Suppose that $A$ is an $\N$ set. Let $\{a_n: n \in \omega\}$ be a
sequence of positive reals witnessing that. We will use the
following notation:
for $n \in \omega$, let $s_n=\sum_{i=0}^n a_i$ and $b_n=a_n/s_n$.

It is well known that $\sum_{n=0}^\infty b_n = \infty$.
Let $\{q_n: n \in \omega\}$ be a sequence of integers such that
$\lim_{n \rightarrow \infty} q_n = \infty$ and
$$\sum_{n=0}^\infty \frac{a_n}{s_n^{1 + \frac{1}{q_n}}} < \infty.$$
Finally let $\varepsilon_n = s_n^{-1/q_n}$.

We will need the following easy lemma:
\begin{lemma}\label{eight}
  Suppose that $Z$ is a finite set of integers, $\varepsilon>0$ and $x
  \in \reals$. Then there exists a set $Z' \subseteq Z$, $|Z'|\geq
  \varepsilon \cdot |Z|$ such that 
$$\forall i,j \in Z' \ \lbv (i-j)x\rbv < 2\varepsilon.$$

\end{lemma}
\Proof Follows immediately from the pigeon-hole principle.~$\QED$

Suppose that $X=\{x_\alpha: \alpha < \delta < \frak t\}$ is a set of
reals.

By induction we will build a sequence $\{\varphi_\alpha : \alpha \leq
\delta \}$ such that for all $\alpha \leq \delta$:
\begin{enumerate}
  \item $ \varphi_\alpha : \dom(\varphi_\alpha)
    \longrightarrow [\omega]^{<\omega}$, $\dom(\varphi_\alpha) \in
    [\omega]^\omega$,
\item $\forall n \in \dom(\varphi_\alpha) \
  \max(\varphi_\alpha(n)) \leq s_n$, 
\item $\forall k \  \lim_{n \in \dom(\varphi_\alpha)}
  |\varphi_\alpha(n)|\cdot \varepsilon_n^k=\infty$,  
\item $\forall \beta>\alpha \ | \dom(\varphi_\beta )\setminus
  \dom(\varphi_\alpha)| < \boldsymbol\aleph_0 $,
\item $ \forall \beta>\alpha \ \forall^\infty n \in \dom(\varphi_\beta)
  \ \varphi_\beta(n) \subseteq \varphi_\alpha(n)$,
\item $ \sum_{n \in \dom(\varphi_\alpha)} b_n =
  \infty$,
\item $ \forall^\infty n \in \dom(\varphi_\alpha) \
  \forall i,j \in \varphi_\alpha(n)\ \lbv (i-j)nx_\alpha \rbv < 
2\varepsilon_n$.
\end{enumerate}

Suppose for a moment that a sequence satisfying these conditions
has been constructed. Let $\varphi = \varphi_\delta$. 
For $n \in \dom(\varphi)$ (sufficiently large) choose two distinct
numbers $i_n,j_n \in 
\varphi(n)$.
Consider the series 
$\sum_{n \in \dom(\varphi)} b_n \lbv (i_n-j_n)nx\rbv$. 
We will show that this series converges for $x \in A \cup X$.

For $x \in A$,
$$\sum_{n \in \dom(\varphi)} b_n \lbv (i_n-j_n)nx\rbv \leq \sum_{n \in \dom(\varphi)}
b_n \cdot 2s_n \cdot \lbv nx\rbv \leq 2 \cdot \sum_{n \in \dom(\varphi)} a_n
\cdot \lbv nx\rbv < \infty.$$

For $x \in X$ we get 
$$\sum_{n \in \dom(\varphi)} b_n \lbv (i_n-j_n)nx\rbv \leq \sum_{n \in \dom(\varphi)}
2b_n \cdot \varepsilon_n = 2\sum_{n \in \dom(\varphi)} \frac{a_n}{s_n^{1 +
    \frac{1}{q_n}}} < \infty.$$ 
Note that by reenumerating terms we can put the series $\sum_{n \in \dom(\varphi)}
b_n \lbv (i_n-j_n)nx\rbv$ in the form required by Definition
\ref{defdef}(1).

Thus to conclude the proof it remains to construct the sequence
$\{\varphi_\alpha : \alpha \leq \delta\}$.

Suppose that $\varphi_\beta$ for $\beta<\alpha$ are given. We will
describe how to find $\varphi_\alpha$.

\vspace{0.1in}

{\sc Case 1} $\alpha$ is limit.
  For $\beta<\alpha $ define $f_\beta  \in \omega^\omega$ by
$$f_\beta (n+1) = 
   \min\left\{k > f_\beta (n): \sum \{b_n : n \in \dom(\varphi_\beta)
   \cap [f_\beta(n),k)\} 
   >1\right\}
   \text{ for } n \in \omega.
$$

   Since $\delta< \frak t \leq \frak b$, let $f \in \omega^\omega$ be
   an increasing function such that 
$$\forall \beta<\alpha   \  \forall^\infty n \ f(n)>f_\beta(n) .
$$

   Let $I_n = \lft1[\widetilde{f}(n), \widetilde{f}(n+1)\rgt1)$ for $n
   \in \omega$, where $\widetilde{f}$ is defined as:
   $\widetilde{f}(0)=f(0)$ and
   $\widetilde{f}(n+1)=f\lft1(\widetilde{f}(n)+1\rgt1)$.

Note that 
$$\forall \beta<\alpha  \ \forall^\infty n \ \widetilde{f}(n) \leq
f_\beta \lft1(\widetilde{f}(n)\rgt1) <  
f_\beta \lft1(\widetilde{f}(n)+1\rgt1) \leq
f\lft1(\widetilde{f}(n)+1\rgt1)=\widetilde{f}(n+1).$$

   It follows
$$\forall  \beta<\alpha   \ \forall^\infty n \ \sum_{j \in I_n \cap
  \dom(\varphi_\beta)} b_j 
   >1.
$$
 
Similarly, we can find a function $g \in \omega^\omega$ such that
\begin{enumerate}
\item $\forall k \ \lim_{n \rightarrow \infty} g(n) \varepsilon^k_n =
  \infty$,
\item $\forall \beta < \alpha \ \forall^\infty n \in \dom(\varphi_\beta) \
|\varphi_\beta(n)| \geq g(n).$
\end{enumerate}

   For $n \in \omega$, let $F \in U_n$ if the 
   following conditions are satisfied:
\begin{enumerate}
\item $F: \dom(F) \longrightarrow [\omega]^{<\omega}$, $\dom(F) \subseteq I_n$,
\item $\max(F(i)) \leq s_i$ for $i \in \dom(F)$,
\item $|F(i)| \geq g(i) \varepsilon_i$ for $i \in \dom(F)$,
\item $\sum_{i \in \dom(F)} b_i > 1/2$,
\item $\forall k \in \dom(F) \ \forall i,j \in F(k) \ \lbv
  (i-j)kx_\alpha \rbv < 2\varepsilon_k$.
\end{enumerate}

   Notice that each set $U_n$ is finite and that $U
   = \bigcup_{n \in  
   \omega} U_n$ can be identified with~$\omega$.

   For $\beta<\alpha $ define
$$X_\beta  = \left\{F  : 
   \exists n \ 
   \lft2( F \in U_n \ \& \ 
     \forall j \in \dom(F) \ F(j) \subseteq \varphi_\beta(j)\rgt2)\right\}.
$$

   It is clear that $|X_\gamma  \setminus X_\beta| <
   \boldsymbol\aleph_0 $
   when $\gamma  \geq \beta$. 
Since $\alpha <\frak t$ we can find a set $X$  such that
   $X \subseteq^\star X_\beta $ for $\beta<\alpha$. We can assume 
   that $X \cap U_n$ consists of at most one point.
   Let
$\varphi'_\alpha = \bigcup X$. To obtain $\varphi_\alpha$ from
$\varphi'_\alpha$ we proceed
exactly as in the successor step below.

\vspace{0.1in}

{\sc Case 2} $\alpha$ is a successor.
To construct $\varphi_\alpha$ from $\varphi_{\alpha-1}$ (or from
$\varphi'_\alpha$ above) we use Lemma \ref{eight}.
In particular, to get $\varphi_\alpha(n)$ we apply Lemma \ref{eight} to
$Z=\varphi_{\alpha-1}(n)$ (or $\varphi'_\alpha(n)$) with
$\varepsilon=\varepsilon_n$ and $x=x_\alpha$. 
Note that 
$$|\varphi_\alpha(n)| \geq |\varphi_{\alpha-1}(n)|  \varepsilon_n
\stackrel{n \rightarrow \infty} \longrightarrow \infty$$
or 
$$|\varphi_\alpha(n)| \geq g(n)  \varepsilon_n
\stackrel{n \rightarrow \infty} \longrightarrow \infty.$$
In either case $\varphi_\alpha$ satisfies the required
conditions.~

This finishes the proof of (1) of Theorem \ref{first}. Part (2) is
proved similarly.
$\QED$
 
\begin{definition}
   Let $\frak t'$ be the least cardinal $\kappa$ such that there exists a 
   family of sequences of positive reals $\{\<a^\alpha_n : n \in 
   \omega\>: \alpha < \kappa\}$ such that
\begin{enumerate}
\item $\sum_{n=0}^\infty a^\alpha_n = \infty$ for all $\alpha$,
\item $\forall \alpha \leq \beta \ 
   \forall^\infty n \ a^\alpha_n \geq \ a^\beta_n$,
\item there is no sequence of positive reals $\<a_n : n \in \omega\>$ 
   such that 
\begin{enumerate}
\item $\forall \alpha  \ 
   \forall^\infty n \ a^\alpha_n \geq a_n$, and
\item $\sum_{n=0}^\infty a_n = \infty$.
\end{enumerate}
\end{enumerate}
\end{definition}

Notice that we proved that:
\begin{corollary}\label{teqtprime}
   $\frak t = \frak t'$.
\end{corollary}

\section{{\normalshape\sf N}$_0$-sets and other sets.}

   To prove the remaining parts of Theorem 1 we need another lemma.
   For an infinite subset $X$ of $\omega$ let $X(n)$ denote the $n$-th
   element of $X$.

\begin{lemma}\label{elias2}
   Let $x_0$ be a real number and let $Z$ be an infinite set of
   natural numbers. There exists an infinite set $X \subseteq Z$
   such that for every infinite set $Y \subseteq X$,
$$\sum_{n=1}^{\infty}\lbv \lft1(Y(n+1)-Y(n)\rgt1)\cdot x_0\rbv<\infty.$$

 \end{lemma}
\Proof 
   Let $x_0$ be a real number and let $Z$ be an arbitrary infinite
   subset of $\omega$. Define by induction a sequence 
   $X_0 \supseteq X_1 \supseteq X_2 \cdots$ such that 
$X_0\subseteq Z$  and 
$$\forall  n \ \forall i,j \in X_n \ \lbv
   |i-j|\cdot x_0\rbv \leq \frac{1}{2^n}.$$
   Suppose that $X_n$ is already defined.
   Define $c:[X_n]^2 \longrightarrow 2$ so that:
\[c(i,j)=\left\{\begin{array}{ll}
   1 & \text{if } \lbv |i-j|\cdot x_0\rbv \leq 2^{-(n+1)}\\
   0 & \text{otherwise}
   \end{array}
   \right. .
\]
   By Ramsey's theorem there exists an infinite set $W$ such that $c$ is 
   constant on $[W]^2$. Note that by Lemma \ref{eight} we have
   $c(i,j)=1$ for all $i,j \in W$. Put $X_{n+1}=W$.

   Let $X \subseteq Z$ be such that $X(n) \in X_{n}$ for all $n$. 
   Then $X$ has the required properties.~$\QED$

   Let $Z$ and $x_0$ be as in the hypotheses of Lemma \ref{elias2}. Set
$${\cal A}_{x_0}=\left\{A\subseteq Z:
\forall B \subseteq A \ \sum_{n=1}^{\infty}\lbv \lft1(B(n+1)-B(n)\rgt1)\cdot
   x_0\rbv<\infty\right\}.$$
   Then Lemma \ref{elias2} asserts that ${\cal A}_{x_0}$ contains a dense
   open subset of ${ P}(Z)/{\text{\sf fin}}$, ordered by almost inclusion.

\begin{lemma}\label{elias3}
   If $A$ is an $\N_0$ set and the infinite set $Z=\<Z(n): n \in \omega\>$
   witnesses that, then for any infinite set $X\subseteq Z$ such that
   $\<X(n+1)-X(n):n\in\omega\>$ is increasing, the sequence
   $\lft1\<X(n+1)-X(n):n<\omega\rgt1\>$, witnesses that $A$ is an $\N_0$-set. 
\end{lemma}

\Proof Note that for every $x \in A$,
$$\lbv \lft1(X(n+1)-X(n)\rgt1)\cdot x \rbv \leq \lbv X(n+1)\cdot x\rbv + 
   \lbv X(n)\cdot x\rbv.~\QED$$

{\sc Proof of theorem \ref{first}(2)-(3)} 
   Let $A$ be an $\N_0$-set and let $Z$ be an infinite set which
   witnesses this. Let $\{x_\alpha: \alpha<\delta<\frak h\}$ be a set
   of real numbers. By Lemma \ref{elias2}, each of the sets ${\cal
     A}_{x_{\alpha}}$ contains a dense open subset of ${
     P}(Z)/{\text{\sf fin}}$. 

   Since $\delta<\frak h$, we find for each $\alpha<\delta$ a set
   $A_{\alpha}\in{\cal A}_{x_{\alpha}}$ and an infinite subset $B$ of $Z$
   such that $B\backslash A_{\alpha}$ is finite for each $\alpha$. By
   further thinning out, we may assume that $\<B(n+1)-B(n):n<\omega\>$
   is strictly increasing. But then, by Lemmas \ref{elias2} and
   \ref{elias3}, the sequence $\<B(n+1)-B(n):n<\omega\>$ witnesses that
   $A\cup\{x_\alpha:\alpha<\delta\}$ is an $\N_0$-set.

   The proof of part (3) is similar.~$\QED$

\section{Arbault and Hardy sets.}
To prove the last two parts  of Theorem \ref{first} we will use the
following lemmas: 
\begin{lemma}[Booth \cite{Booth74}]\label{booth}
Suppose that $X$ is a set of reals of size $< \frak s$. Let $f_n: X
\longrightarrow [0,1]$ 
for $n \in \omega$ be a sequence of functions. Then there exists a
sequence $\<n_k : k \in \omega\>$ 
such that $\<f_{n_k}: k \in \omega\>$ converges pointwise on $X$.
\end{lemma}
\Proof
For $q \in \rationals \cap [0,1]$ and $x \in X$ let
$$A^{\mathord{<}}_{q,x} = \{n : f_n(x)\leq q\} \text{ and } 
A^{\mathord{>}}_{q,x} = \{n : f_n(x)> q\}.$$
Let $Y=\<n_k: k \in \omega\}$ be a set witnessing that 
$\{A^{\mathord{<}}_{q,x}, A^{\mathord{>}}_{q,x}: q \in \rationals \cap
[0,1], x \in X\}$ 
is not a splitting family. 
It is easy to see that this is a sequence of the sort we are looking
for.~$\QED$ 

We have the following analog of Lemma \ref{elias3}:
\begin{lemma}\label{e6}
   If $A$ is an Arbault set and the infinite set $Z=\<Z(n): n \in \omega\>$
   witnesses that, then for any infinite set $Y\subseteq Z$ such that
   $\lft1\<Y(n+1)-Y(n):n\in\omega\rgt1\>$ is increasing, the sequence
   $\lft1\<Y(n+1)-Y(n):n<\omega\rgt1\>$, witnesses that $A$ is an
   Arbault set. 
\end{lemma}
\Proof Use the fact that:
$$\lft1|\sin(x-y)\rgt1| = \lft1| \sin(x)\cos(y) - \sin(y)\cos(x)\rgt1|
\leq
\lft1|\sin(x)\rgt2|+ \lft1| \sin(y)\rgt1|.~\QED$$

{\sc Proof of theorem \ref{first}(4)-(5)} Suppose that $A$ is an Arbault 
set and
$|X| < \frak s$. Assume that sequence $Z = \<n_k : k \in \omega\>$
witnesses that $A$ is an Arbault set.

Let $f_k(x) = \sin(\pi \cdot n_k x)$ for $x \in X$.
By Lemma \ref{booth} there exists a set $Y \subseteq Z$ such that
the sequence $\lft2\< \sin\lft1(\pi \cdot Y(n) \cdot
x\rgt1)\rgt2\>$ 
converges for every $x \in X$.

In particular, $\<Y(n+1)-Y(n): n \in \omega\>$ witnesses that $A \cup
X$ is an Arbault set.

\vspace{0.1in}
Theorem \ref{first}(5) follows immediately from the fact 
 that Arbault sets are $\H_\sigma$ sets. In other words, a set of size
 $< \frak s$ is a Arbault set and so a $\H_\sigma$ set. The union of two
 $\H_\sigma$ sets is an $\H_\sigma$ set.~$\QED$

\ifx\undefined\bysame
\newcommand{\bysame}{\leavevmode\hbox to3em{\hrulefill}\,}
\fi


\begin{thebibliography}{1}

\bibitem{bary64}
N.K. Bary, {\em A treatise on trigonometric series}, MacMillan Book Co., New
  York (1964), Translated by M.F. Mullins.

\bibitem{Booth74}
D.~Booth, {\em A {B}oolean view of sequential compactness}, Fundamenta
  Mathematicae {\bf 85} (1974), 99--102.

\bibitem{BuBu}
Z.~Bukovsk\'a and L.~Bukovsk\'y, {\em Adding small sets to an ${N}$-set} (1993),
  preprint.

\bibitem{Buki}
Lev Bukovsk\'y, {\em Thin sets related to trigonometric series}, Israel
  {M}athematical {C}onference {P}roceedings, Bar Ilan University (1992),
  pp.~107--119.

\bibitem{elias93}
Peter Eli\'a\v{s}, {\em Systematizacia tried malych mnozin suvisicich s
  absolutnou konvergenciou trigonometrickych radov}, 1993, diplomov\'a pr\'aca,
  University P.J.\v Saf\'arik, Ko\v{s}ice.

\bibitem{Kch}
N.N. Kholshchevnikova, {\em On the properties of thin sets for trigonemetric
  series and certain other series}, Dokl. Russian Acad. Sci. (1992), no.~4--6,
  446--449.

\bibitem{vdo:inttop}
E.~K. van Douwen, {\em The integers and topology}, Handbook of Set Theoretic
  Topology (Amsterdam) (K.~Kunen and J.~E. Vaughan, eds.), North-Holland,
  Amsterdam(1984), pp.~111--167.

\bibitem{Vau90Sma}
Jerry E. Vaughan, {\em Small uncountable cardinals and topology}, {O}pen problems
  in topology (Amsterdam), North Holland, Amsterdam (1990).

\end{thebibliography}
\end{document}